\documentclass[11pt]{amsart}

\usepackage[utf8]{inputenc}
\usepackage[T1]{fontenc}
\usepackage[english]{babel}

\usepackage{amsmath}
\usepackage{amssymb}
\usepackage{amsthm}
\usepackage{amsfonts}
\usepackage{amscd}

\usepackage{graphicx}
\usepackage{multirow}
\usepackage{booktabs}
\usepackage{array}

\usepackage{enumerate}
\usepackage{url}
\usepackage[colorlinks=true, linkcolor=blue, citecolor=blue, urlcolor=blue]{hyperref}
\usepackage{tikz}

\theoremstyle{plain}
\newtheorem{theorem}{Theorem}[section]

\theoremstyle{definition}
\newtheorem{definition}[theorem]{Definition}
\newtheorem{example}[theorem]{Example}

\theoremstyle{remark}

\title[Symmetry in Tree Parking]{Symmetry in Tree Parking Distributions}

\author{Amanuel T. Getachew}

\subjclass[2020]{Primary 05A15; Secondary 05A19, 05C05}
\keywords{Fuss-Catalan numbers, symmetric functions, caterpillar trees, parking functions}

\date{\today}

\begin{document}

\begin{abstract}
In this paper, we explore parking distributions on caterpillar trees, focusing on two primary statistics: the number of lucky cars and the frequency with which cars prefer specific parking spaces. We use first-return decomposition to reveal a symmetry in their joint distribution and develop a $q,t$-analog of the Fuss-Catalan generating function. We prove that this generating function exhibits specific symmetry and satisfies a functional equation. Additionally, we extend our findings to any $m$ statistics that satisfy certain criteria, presenting a concrete example of such $m$ statistics to illustrate the broader applicability of our results.
\end{abstract}

\maketitle

\section{Introduction}

\label{intro}
Parking functions are a fundamental concept in combinatorics, with wide-ranging applications and characterizations in various areas of mathematics and computer science. Several generalizations and refinements of parking functions have been reported in the literature, one of which involves non-decreasing parking functions on rooted trees, also known as tree parking distributions. In this paper, we explore two commonly studied statistics in tree parking distributions and their symmetry, which leads to a breakdown of the generating functions that enumerate these combinatorial objects.

\subsection{Parking Functions, Parking Distributions and Caterpillar Trees}
\label{subsect1.1}

Parking functions were introduced by Konheim and Weiss to study hashing protocols in \cite{Konheim_1966}. Consider a parking lot with $m$ parking spaces and $n$ drivers. The drivers come sequentially and search for unoccupied parking spaces. The $i$th driver has a designated preferred parking space $p_i$. If $p_i$ is unoccupied upon arrival, the driver parks there. However, if $p_i$ is already occupied, the driver sequentially checks the next available spaces, $pi+1,pi+2$, and so forth, until an unoccupied space is found, where they then park. If the driver reaches the last parking space and it is occupied, the driver will leave the parking lot. We call this process the parking process.

A parking function of length $n$ is a sequence of $n$ parking preferences $(p_1, p_2, \dots, p_n)$ such that all $n$ drivers are able to park. There are \((n + 1 - m) \cdot (n + 1)^{(m-1)}\) parking functions, for \(n\) parking spaces and \(0 \leq  m \leq n\) drivers.

Since their introduction, parking functions have been 
studied extensively and connections to various other combinatorial objects such as hyperplane arrangements 
\cite{Stanley1998}, acyclic functions \cite{Franon1975}  and noncrossing partitions \cite{StanNonCrossing} have been revealed. Moreover, the notion of parking functions has been generalized in several ways, yielding, e.g., \((a, b)\)-parking functions \cite{YanGeneral} and  graph-parking functions \cite{Gpark}.

A parking distribution is a nondecreasing parking function. Mathematically, a parking distribution $\mathbf{a} = (a_1 , . . . , a_n)$ is an non-decreasing sequence of positive integers in $[n]$ that satisfies the inequality

\[a_{i} \leq i,\]
or equivalently, 
\begin{equation}
    \sum_{1\leq i \leq k} |\{j:a_j = i \}| \geq k \ \text{for all} \  1\leq k \leq n
    \label{defeq1}
\end{equation}
with the equality holding for $k = n.$

It is  obvious  that a parking function $(p_1, p_2, \dots, p_n)$ can be uniquely determined by a parking distribution - permutation pair, ($\mathbf{a}, \sigma)$  as $p_i = a_{\sigma(i)}$. Moreover, parking distributions,  as it is  shown in the subsequent sections, are enumerated by Catalan numbers which makes them useful for developing theories on other combinatorial structures. For example,  \cite{hopfalg} utilized a bijection between noncrossing partition and parking distributions to develop the Hopf algebra of parking functions.

Let $T$  be called a sink tree if $T$ is directed tree and if there is a path from every node   or vertex in $T$ to the root (called sink) of $T$. In \cite{Butler_2017}, the authors generalized the parking process described above to sink trees, i.e., each car starts from the node of its preference and  moves along the tree edge until it finds an available node or exits through the sink. A nondecreasing sequence of parking preferences $\mathbf{a} = (a_1, \dots, a_n)$  is called $T-$parking distribution if all $n$ cars can park on $T$. 

When $T$ is a path graph, $\mathbf{a}$ is an ordinary parking distribution defined by equation \ref{defeq1}. An interesting case is when  $T$ is a tree formed by zero or more nodes directly connected to (possibly the middle of) a  path graph called the \textit{backbone}. If $T$ is formed by connecting $b_i-1$ nodes to the $i$-th node of the backbone, then the set of parking distributions on  $T$ is in bijection with lattice paths strictly on the left of the boundary $\{(0, b_1), (1, b_1 + b_2), \dots \}$ \cite{Butler_2017}. In this paper, we study the case when $b_1 = m$ and $b_i = m-1$ for $i > 1$ and some constant $m$, i.e., when $T$ is caterpillar tree. As we shall see, the number of parking distributions on these trees, is enumerated by the $m$-Fuss-Catalan numbers which are of major interest in enumerative combinatorics \cite{Multivar,density}.
\subsection{Fuss-Catalan numbers}
\label{subsect1.2}
The Fuss-Catalan (also called $m$-Catalan or $m$-Fuss-Catalan) numbers are the numbers of the form
\begin{equation}
  C^{(m)}_n = \frac{1}{mn + 1} \binom{mn + n}{n}
  \label{fuss-catalan}
\end{equation}
for $n \geq 0$ and for a fixed positive integer $m$. Here are some of the combinatorial interpretations of the $m$-Fuss-Catalan numbers:
\begin{enumerate}
  \item The number of lattice paths from $(0,0)$ to $(n, mn)$ using steps $(1,0)$ and $(0,1)$ that do not go above the line $y = mx$.
  
  \item The number of ways to parenthesize a product of $mn+1$ factors using only $m$-ary operations (operations taking $m$ operands).
  
  \item The number of full $m$-ary trees with $n$ internal nodes, where each internal node has exactly $m$ children.
\end{enumerate}

Their generalizations have been extensively studied in literature \cite{Multivar,hypergraph,density}.  It is relatively well-known \cite{density} that the $m$-Fuss-Catalan numbers are enumerated by the generating function 
$$\mathcal{B}_{m}(x) = \sum_{n \geq 0} C^{(m)}_n x^n  $$
which also satisfies the functional equation

\begin{equation}
  \mathcal{B}_{m}(x) = 1 + x\mathcal{B}_{m}^{m+1}(x).
  \label{funceq}
\end{equation}

\newcommand{\ip}{\mathbf{p}}
\subsection{$\mathrm{luck}$ and $\omega_j$} In the classical parking process, a driver is called \textit{lucky} if it gets to park in its preferred spot. This statistic is of interest in the literature (see \cite{luckycars}, for example). By extension, in caterpillar parking distributions, we call a car lucky if it \textbf{prefers to park at a node in the backbone} and gets to park there. For a caterpillar parking distribution $\ip$, we denote the number of lucky drivers by $\mathrm{luck}(\ip)$. Another common statistic in enumeration is the frequency statistic, $\omega_j$. We define $\omega_k(\ip)$ as the number of drivers that prefer to park at space $k$.

In \cite{Garsia1996ARQ}, the authors introduced a family of functions, $C_n(q, t)$, which were later proven to be polynomials in \cite{Haiman_2002}, called the $q, t$-Catalan sequences. It was proven in \cite{Garsia_2001} that
$$C_n(q, t) = \sum_{D \in \mathcal{D}_n} q^{\mathrm{dinv}(D)}t^{\mathrm{area}(D)},$$
where $\mathcal{D}_n$ is the set of Dyck paths of semi-length $n$, $\mathrm{dinv}(D)$ is the number of diagonal inversions (pairs of cells $(i,j)$ and $(i',j')$ with $i < i'$, $j > j'$, and both cells under the path), and $\mathrm{area}(D)$ is the number of unit squares between the path and the $x$-axis. It also follows from the definition in \cite{Garsia1996ARQ} that $C_n(q, t) = C_n(t, q)$. Motivated by these results and the connection between Dyck paths and parking distributions (as we shall see later), in this paper, we introduce another $q, t$-Catalan analog for parking distributions and extend the result to $m > 2$. To do so, we define $\gamma_n$ as
$$\gamma^{(m)}_n(q, t) = \sum_{\ip} q^{\mathrm{luck}(\ip)}t^{\omega_1(\ip)},$$
where the sum extends over all tree parking distributions on an $m$-regular caterpillar tree (defined later). We prove that this $q, t$-analog satisfies that

\begin{enumerate}
  \item  $\gamma^{(m)}_n(1, 1) = C_n^{(m)}$
  \item  $\gamma^{(m)}_n(q, t) = \gamma^{(m)}_n(t, q)$
  \item  $\sum_{n \geq 0} \gamma^{(m)}_n(q, 1) x^n = \frac{1}{1 - qx\mathcal{B}_m(x)}$
  \item $\gamma^{(m)}_n(q, t)$ is a linear combination of complete homogeneous polynomials, $h_i(q, t)$ for $0 \leq i \leq n$. 
  \end{enumerate}
  
  This will be the first main result of our paper. To prove this result, we introduce the notion of fixed points of a parking distribution of various types and extend the notion of first-return decomposition of Dyck paths introduced in \cite{Deutsch_1999} to tree parking distributions. We show that the $q, t$-analogue of the generating functions $\mathcal{B}_m(x)$, $\mathcal{B}_m(x; q, t)$ satisfies the functional equation
  $$\mathcal{B}_m(x; q, t) = 1 + xqt\mathcal{B}_m(x;q, 1)\mathcal{B}_m(x;1, t)\mathcal{B}_m^{m-1}(x),$$
  which proves our main result.
\subsection{More statistics} 
We generalize our results for $\mathrm{luck}$ and $\omega_1$ to $m+1$ hypothetical statistics, $\mathcal{S}_0, \dots, \mathcal{S}_{m}$, that satisfy certain criteria. Consequently, we show that the symmetric joint distribution also applies to these statistics. We then prove that $\mathrm{luck}$ and $\omega_k$ for $1 \leq k \leq m$ satisfy the criteria we will define, obtaining a concrete example for the generalization. 

The rest of the paper is organized as follows. In section \ref{secParD}, we provide formal definitions for caterpillar trees and parking distributions and provide examples of them. In section \ref{secUpar}, we give a simple bijection linking caterpillar tree parking distributions and $\mathbf{u}$-parking distributions. In section \ref{secFTD}, we introduce the first-return decomposition, which is our main tool for the main theorem. In section \ref{secEqui} and \ref{secSymm}, we prove our first main result that $\mathrm{luck}$ (lucky drivers) and $\omega_1$ (first-space-loving drivers) have a symmetric joint distribution. Our second main result, a generalization, will be presented in section \ref{secGen}. Discussion and unsolved problems are addressed in the concluding section.

\section{Parking Distributions on Caterpillar Trees \label{secParD}}
We begin by formally defining parking distributions:
\begin{definition}
  \label{def11}
    Let $T$ be a rooted tree with vertices labeled $1, 2, \dots , n$ and let $T_k$ be the sub-tree of $T$ rooted at $k$. Let $p_T: V(T) \to \mathbb{N}$ be a distribution of indistinguishable balls on the vertices of $T$. We say that $p_T$ is a parking distribution if
    $$\sum_{i \in T_k}p_T(i) \geq |T_k|$$ 
    for all vertices $k$.
\end{definition}

As mentioned above, one can easily observe that this definition is the same as the definition of classical parking distributions given by equation \ref{defeq1} when $T$ is a path.

The authors of \cite{Butler_2017} also discussed parking distributions on a special type of digraphs called \textit{caterpillar trees}. Here, we define a subclass of those.

\newcommand{\Cat}[2]{\mathrm{Cat}_{#1}(#2)}

\begin{definition}
  \label{def12}    
  A tree $T$ is called an \textbf{$m$-caterpillar} (or \textbf{$m$-regular caterpillar}) tree of length $n$ for $m \geq 1$, if there is a (directed) path $v_1 \to v_2 \to \dots \to v_n$ such that $v_2$ is connected to $m$ leaves including $v_1$ and each $v_i$ for $i \geq 3$ is connected to $m-1$ leaves. 
\end{definition}
We say the vertices $v_1, v_2, \dots, v_n$ form the backbone of the caterpillar tree. Let $d(v)$ denote the depth of a vertex $v \in V(T)$ and $\deg v$ denote the in-degree of $v$. There exists a partial order, $<_*$, of the vertices of $T$, namely for two vertices $u, v \in V(T)$, $u <_* v$ if $d(v) < d(u)$ if $d(u) \neq d(v)$ or $\deg v < \deg u$ if $d(v) = d(u)$. We refer to this partial order to define labelling of nodes in $T$ as follows. All nodes are labelled $1$ to $mn- m + 1$. If $<_*$ is defined for two nodes  $u, v$, the label of node $u$ shall be less than that of $v$ if and only if $u <_* v$. Otherwise, $u$ and $v$ shall have distinct labels in any order.  We denote an $m$-regular caterpillar tree labelled this way with $\Cat{m}{n}$ and the set of all parking distributions on $\Cat{m}{n}$ with $\mathcal{PK}_m(n)$. In the remaining part of the paper, we focus on parking distributions in which the drivers occupy the parking spaces of $\Cat{m}{n}$ in this order.

    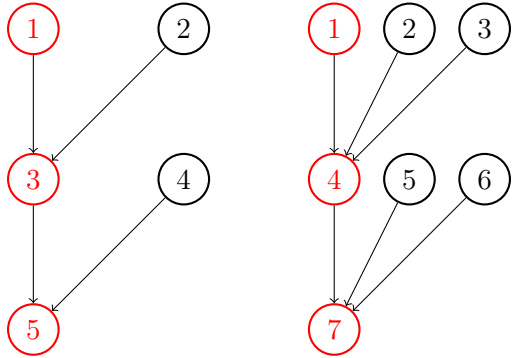
\begin{figure}
  \centering
   \caption{$\Cat{2}{3}$ and $\Cat{3}{3}$; The red vertices form the backbones.}
   \vspace{.5cm}
  \begin{tikzpicture}[every node/.style={circle, thick, draw}]
      \node[red] (A) at (0,0) {5};
      \node[red] (B) at (0,2) {3};
      \node (C) at (2,2) {4};
      \node[red] (D) at (0,4) {1};
      \node (E) at (2, 4) {2};
      
      \draw [->] (B) -- (A);
      \draw [->] (C) -- (A);
      \draw [->] (D) -- (B);
      \draw [->] (E) -- (B);
  
       \node[red] (A7) at (4,0) {7};
      \node[red] (A4) at (4,2) {4};
      \node (A5) at (5,2) {5};
      \node (A6) at (6,2) {6};
      \node[red] (A1) at (4,4) {1};
      \node (A2) at (5,4) {2};
      \node (A3) at (6,4) {3};
      
      \draw [->] (A4) -- (A7);
      \draw [->] (A5) -- (A7);
      \draw [->] (A6) -- (A7);
      \draw [->] (A1) -- (A4);
      \draw [->] (A2) -- (A4);
      \draw [->] (A3) -- (A4);
    
  \end{tikzpicture}
 
  \end{figure}

\begin{example}
    The 12 parking distributions on $\Cat{2}{3}$ are shown on the next page.
\begin{table}[h]
  
    \centering
    \caption{Parking Distributions on $\Cat{2}{3}$}
    \label{tab:pf}
    \begin{tabular}{cc}    
        $(1, 1, 1,2,4)$ & $ (1,2,2,3,4)$\\
        $(1, 1, 2,2,4)$ & $(1,2,3,4,4)$\\
        $(1, 1, 2, 3, 4)$ & $(1,2,2,4,5)$\\
        $(1, 1, 2, 4,4)$ & $(1,2,3,3,4)$\\
        $(1, 1, 2,4,5)$ & $(1,2,3,4,4)$\\
        $(1, 2, 2, 2, 4)$ & $(1,2,3,4,5)$\\

    \end{tabular}
\end{table}
    
\end{example}

\newcommand{\luck}{\mathrm{luck}}
\newcommand{\freq}{\omega}
\newcommand{\p}{\mathbf{p}}
\newcommand{\PK}[2]{\mathcal{PK}_{#1}(#2)}

It is well-known that classical parking distributions are enumerated by Catalan numbers. This can be shown via bijection with Dyck paths of semi-length $n$ defined as follows: For $\p \in \PK{1}{n}$, the word $N^{\freq_1(\p)}E\cdots N^{\freq_n(\p)}E$ is a Dyck word of semi-length $n$.

In \cite{Butler_2017}, the authors show that for a sequence $(a_2, \dots, a_n)$ and a caterpillar tree $T$ with $a_i-1$ leaves connected to the $i$-th backbone vertex $v_i$ for $i > 2$, and $a_2$ leaves to $v_2$, the number of parking distributions on $T$ is given by
$$\mathrm{det}\bigg[ {{\sum_{j \leq r} a_{r+1}} \choose {s - r + 1}} \bigg]_{1 \leq r, s \leq n-1}.$$
When $a_2 = m$ and $a_i = m - 1$ for $i > 2$, we have
$$|\PK{m}{n}| = \dfrac{1}{mn - m + 1} {{mn} \choose {n}},$$

We study caterpillar parking distribution by defining a bijection with similar combinatorial structures that are enumerated by the $m$-Fuss-Catalan numbers described below.
\newcommand{\uu}{\mathbf{u}}
\section{$\mathbf{u}$-Parking Distributions}\label{secUpar}
Another generalizations of classical parking functions are the $\uu$-parking functions \cite{Kung2003}. We define the parking distribution analoges, $\uu$-parking distributions.

\begin{definition}
  Let $\mathbf{u} = (u_1, u_2, \dots)$ be an increasing sequence of positive integers. An increasing sequence of positive integers $\mathbf{a} = (a_1, \dots, a_n)$ such that $1 \leq a_1 \leq u_1$ and  $a_{i-1} \leq a_i \leq u_i$ for $i > 1$ is called a $\uu$-parking distribution. We denote the set of all $\uu$-parking distributions of length $n$ by $\PK{}{n;\uu}.$
\end{definition}
When $\uu = (1, 2, \dots)$, we have $\PK{}{n; \uu} = \PK{1}{n}.$ 

\newtheorem{prop}{Proposition}
\begin{prop}
      Let $\uu = (1, m+1, 2m + 1, \dots)$. Then we have
$$|\PK{}{n;\uu}| =  \dfrac{1}{mn - m + 1} {{mn} \choose {n}}$$
\end{prop}

\begin{proof}



We show that $\PK{}{n;\uu}$ is in bijection with the set of lattice paths from $(0, 0)$ to $(mn-m, mn - m)$ with north-steps $N = (0, m)$ and east-steps $E = (1, 0)$ and staying above the line $y = x$ which are enumerated by $C_{n-1}^{(m+1)}$ as stated in $\ref{subsect1.2}$. Let $\p = (p_1, \dots, p_n) \in \PK{}{n;\uu}$. Then the map $\p \mapsto NE^{p_2 - p_1} \cdots NE^{p_n - p_{n-1}}$ is a bijection between the two sets where $E^k$ represents $k$ consequtive $E$-steps.

  \end{proof}

  \begin{prop}
    \label{prop:iso}
    Let $\uu = (1, m+1, 2m+1, \dots)$. The sets $\PK{}{n ; \uu}$ and $\PK{m}{n}$ are isomorphic.
    \end{prop}
    \begin{proof}
    We define a bijection $\theta$ from $\PK{}{n ; \uu}$ to $\PK{m}{n}$ as follows: if $\p = (p_1, \dots, p_n) \in \PK{}{n ; \uu}$, then let $\theta(\p)$ be an increasing sequence formed by augmenting the numbers in $[mn - m + 1] - \{mi + 1: i \mod n \neq 0 \}$ to $\p$ and then sorting it in increasing order.
    \end{proof}
    \begin{table}[h]
      
      \centering
      \caption{$\theta$}
      \label{tab:theta}
      
      \begin{tabular}{c|c}    
          $\p \in \PK{}{3;(1, 3, \dots)}$ & $\theta(\p) \in \PK{2}{3}$ \\
          \hline
          \hline
          $(1, 1, 1)$ & $(1, 1, 1, 2, 4) $\\
          $(1, 1, 2)$ & $(1, 1, 2, 2, 4) $\\
          $(1, 1, 3)$ & $(1, 1, 2, 3, 4) $\\
          $(1, 1, 4)$ & $(1, 1, 2, 4, 4) $\\
          $(1, 1, 5)$ & $(1, 1, 2, 4, 5) $\\
          $(1, 2, 2)$ & $(1, 2, 2, 2, 4) $\\
          $(1, 2, 3)$ & $(1, 2, 2, 3, 4) $\\
          $(1, 2, 4)$ & $(1, 2, 2, 4, 4) $\\
          $(1, 2, 5)$ & $(1, 2, 2, 4, 5) $\\
          $(1, 3, 3)$ & $(1, 2, 3, 3, 4) $\\
          $(1, 3, 4)$ & $(1, 2, 3, 4, 4) $\\
          $(1, 3, 5)$ & $(1, 2, 3, 4, 5) $\\
  
      \end{tabular}
      
  \end{table}

  Let $h_{n, k, r}^{(m)} = |\PK{}{n; \uu}|$ where $u_i = m(i+k-1) - r$. Let $[x^n] H_m(x; k,r) = h_{n, k, r}^{(m)}.$ Then we have the following theorem:
  \begin{theorem}
    \label{thm:rec}
    $$H_m(x;k,r) = \mathcal{B}_m^{mk - r}(x)$$
    \end{theorem}

\begin{proof}
      We induct on $k$ and $r$. Let $\mathbf{p}^{(j)} = (p_1^{(j)}, \dots, p_n^{(j)})$ be the number of increasing sequences of positive integers such that $p_s^{(j)} \leq m(s+k-1) - (r+1)$ for $1 \leq s \leq j-1$ and $p_j^{(j)} = m(j+k-1) - r$. Since $\mathbf{p}^{(j)}$ is increasing, the subsequence $(p_j^{(j)}, p_{j+1}^{(j)}, \dots, p_n^{(j)})$ should then satisfy the condition

    $$p_{j+l-1}^{(j)} \leq p_{j+l}^{(j)} \leq m(j + l + k - 1) - r$$
    $$\iff p_{j+l-1}^{(j)} - m(j + k - 1) - r + 1 \leq p_{j+l}^{(j)} - m(j + k - 1) - r + 1 \leq ml + 1$$
    for $1 \leq l \leq n - j$. Since $p_j^{(j)} = m(j + k - 1) - r$, this subsequence is enumerated by $h_{n-j+1, 1, m-1}^{(m)}$. 

    On the other hand, if $0 \leq r < m-1$, the increasing subsequence to the left of $j$, $(p_1^{(j)}, \dots, p_{j-1}^{(j)})$, belongs to $\PK{}{j-1; \mathbf{u}'}$, where $\mathbf{u}' = (u_1', u_2', \dots)$ and $u_s' \leq m(s + k - 1) - (r + 1)$, which are enumerated by $h_{j-1, k, r+1}^{(m)}$. Hence, for $0 \leq r < m - 1$, the total number of possible sequences $\mathbf{p}^{(j)}$ is $h_{j-1, k, r+1}^{(m)} \cdot h_{n-j+1, 1, m-1}^{(m)}$. Since $j$ can take any value in $[n+1]$, we have

    $$h_{n, k, r}^{(m)} = \sum_{j = 1}^{n+1} h_{j-1, k, r+1}^{(m)} \cdot h_{n-j+1, 1, m-1}^{(m)} = \sum_{j = 0}^{n} h_{j, k, r+1}^{(m)} \cdot h_{n-j, 1, m-1}^{(m)}.$$

    If $r = m-1$, the subsequence to the left of $j$, $(p_1^{(j)}, \dots, p_{j-1}^{(j)})$, belongs to $\PK{}{j-1; \mathbf{u}''}$, where $\mathbf{u}'' = (u_1'', u_2'', \dots)$ and $u_s'' \leq m(s + k - 1) - m = m(s + k - 1 - 1)$, which is enumerated by $h_{j-1, k-1, 0}^{(m)}$. By a similar reasoning, in this case, we have

    $$h_{n, k, m-1}^{(m)} = \sum_{j = 0}^{n} h_{j, k-1, 0}^{(m)} \cdot h_{n-j, 1, m-1}^{(m)}.$$
    Note that $h_{n, 1, m-1}^{(m)}$ is the $n$-th Fuss-Catalan number $C_n^{(m)}$. Applying Cauchy's product rule, we have

    $$H_m(x; k, r) = \begin{cases}
        H_m(x; k, r+1) \cdot \mathcal{B}_m(x) & \text{if } 0 \leq r < m-1, \\
        H_m(x; k-1, 0) \cdot \mathcal{B}_m(x) & \text{if } r = m-1
    \end{cases}$$ which proves the theorem.
\end{proof}

\section{First-Return Decomposition}\label{secFTD}

  The general idea of the first-return decomposition of combinatorial structures enumerated by the Catalan numbers has been used in the literature (see \cite{flattendcatalan,Deutsch_1999} for example). We extend this idea to $(1, m+1, 2m+1, \dots)$-parking distributions.

  \begin{definition}
    Let $\uu = (1, m+1, 2m+1, \dots)$ and let $\p = (p_1, \dots, p_n) \in \PK{}{n; \uu}$. We define the first fixed point of $\p$ of type $\ell$ as the smallest $k > 1$ such that $m(k - 2) + 1 + \ell \leq p_k \leq m(k-1) + 1$, where $1 \leq \ell \leq m$. The first fixed point of $\ip$ of type $m$ is the smallest $k > 1$ such that $p_k = k$.
\end{definition}

\begin{definition}
Let $i_1 \leq i_2 \leq \dots \leq i_m$ be of the first fixed points of type $1, 2, \dots, m$, respectively.  If $\delta_k$ is the shift-by-$k$ operator, that is, $\delta_k(a_1, a_2, \dots) = (a_1 - k, a_2 - k, \dots)$, then the decomposition of $\p$ into $m+1$ (possibly empty) $\uu$-parking distributions as $(\p_1, \p_2, \dots, \p_{m+1})$ such that
  $$
  \begin{array}{rcl}
      \p_1 & = & \delta_{p_2 - 1}(p_2, \dots, p_{i_1 - 1}) \\
      \p_2 & = & \delta_{p_{i_1} - 1}(p_{i_1}, \dots, p_{i_2 - 1}) \\
      & \vdots & \\
      \p_m & = & \delta_{p_{i_{m-1}} - 1}(p_{i_{m-1}}, \dots, p_{i_m - 1}) \\
      \p_{m+1} & = & \delta_{p_{i_m} - 1}(p_{i_m}, \dots, p_n) \\
  \end{array}
  $$
  is called the first-return decomposition of $\p$.
\end{definition}

For example, if $m = 3$ and $\p = (1, 2, 5, 10, 10, 16)$. $i_1 = 2$,  $i_2= i_3 = 4$. $\p_1 = \varepsilon$, $\p_2 = (1,4)$, $\p_3 = \varepsilon$ and $\p_4 = (1, 1, 7)$. 
We list the first-return decompositions of \(\PK{}{3;(1,3,\dots)}\) and \(\PK{}{3; (1, 4, \dots)}\) in tables \ref{tab:tbldec} and \ref{tab:tbldec2}.

\begin{table}[h]
  
  \centering
  
  \caption{First-Return Decomposition of $\PK{}{3;(1,3,\dots)}$}
  \label{tab:tbldec}
  \begin{tabular}{c|ccc}    
    
      $\p$ & $\p_1$ & $\p_2$ & $\p_3$ \\
      \hline
      \hline
      $(1, 1, 1)$& $ (1, 1)$& $ \varepsilon$& $ \varepsilon$\\
      $(1, 1, 2)$& $ (1, 2)$& $ \varepsilon$& $ \varepsilon$\\
      $(1, 1, 3)$& $ (1, 3)$& $ \varepsilon$& $ \varepsilon$\\
      $(1, 1, 4)$& $ (1)$& $ (1)$& $ \varepsilon$\\
      $(1, 1, 5)$& $ (1)$& $ \varepsilon$& $ (1)$\\
      $(1, 2, 2)$& $ \varepsilon$& $ (1, 1)$& $ \varepsilon$\\
      $(1, 2, 3)$& $ \varepsilon$& $ (1, 2)$& $ \varepsilon$\\
      $(1, 2, 4)$& $ \varepsilon$& $ (1, 3)$& $ \varepsilon$\\
      $(1, 2, 5)$& $ \varepsilon$& $ (1)$& $ (1)$\\
      $(1, 3, 3)$& $ \varepsilon$& $ \varepsilon$& $ (1, 1)$\\
      $(1, 3, 4)$& $ \varepsilon$& $ \varepsilon$& $ (1, 2)$\\
      $(1, 3, 5)$& $ \varepsilon$& $ \varepsilon$& $ (1, 3)$\\
      \vspace*{0.5cm}

  \end{tabular}

\end{table}

\begin{table}[]
  \centering
  \caption{First-Return Decomposition of $\PK{}{3;(1,4,\dots)}$}
  \label{tab:tbldec2}
   \begin{tabular}{c|cccc}    
      $\p$ & $\p_1$ & $\p_2$ & $\p_3$ & $\p_4$ \\
      \hline
      \hline
      $(1, 1, 1)$ & $(1, 1)$ & $ \varepsilon$ & $\varepsilon$ & $\varepsilon$ \\
      $(1, 1, 2)$ & $(1, 2)$ & $ \varepsilon$ & $\varepsilon$ & $\varepsilon$ \\
      $(1, 1, 3)$ & $(1, 3)$ & $ \varepsilon$ & $\varepsilon$ & $\varepsilon$ \\
      $(1, 1, 4)$ & $(1, 4)$ & $ \varepsilon$ & $\varepsilon$ & $\varepsilon$ \\
      $(1, 1, 5)$ & $(1)$ & $ (1)$ & $\varepsilon$ & $\varepsilon$ \\
      $(1, 1, 6)$ & $(1)$ & $ \varepsilon$ & $(1)$ & $\varepsilon$ \\
      $(1, 1, 7)$ & $(1)$ & $ \varepsilon$ & $\varepsilon$ & $(1)$ \\
      $(1, 2, 2)$ & $\varepsilon$ & $ (1, 1)$ & $\varepsilon$ & $\varepsilon$ \\
      $(1, 2, 3)$ & $\varepsilon$ & $ (1, 2)$ & $\varepsilon$ & $\varepsilon$ \\
      $(1, 2, 4)$ & $\varepsilon$ & $ (1, 3)$ & $\varepsilon$ & $\varepsilon$ \\
      $(1, 2, 5)$ & $\varepsilon$ & $ (1, 4)$ & $\varepsilon$ & $\varepsilon$ \\
      $(1, 2, 6)$ & $\varepsilon$ & $ (1)$ & $(1)$ & $\varepsilon$ \\
      $(1, 2, 7)$ & $\varepsilon$ & $ (1)$ & $\varepsilon$ & $(1)$ \\
      $(1, 3, 3)$ & $\varepsilon$ & $ \varepsilon$ & $(1, 1)$ & $\varepsilon$ \\
      $(1, 3, 4)$ & $\varepsilon$ & $ \varepsilon$ & $(1, 2)$ & $\varepsilon$ \\
      $(1, 3, 5)$ & $\varepsilon$ & $ \varepsilon$ & $(1, 3)$ & $\varepsilon$ \\
      $(1, 3, 6)$ & $\varepsilon$ & $ \varepsilon$ & $(1, 4)$ & $\varepsilon$ \\
      $(1, 3, 7)$ & $\varepsilon$ & $ \varepsilon$ & $(1)$ & $(1)$ \\
      $(1, 4, 4)$ & $\varepsilon$ & $ \varepsilon$ & $\varepsilon$ & $(1, 1)$ \\
      $(1, 4, 5)$ & $\varepsilon$ & $ \varepsilon$ & $\varepsilon$ & $(1, 2)$ \\
      $(1, 4, 6)$ & $\varepsilon$ & $ \varepsilon$ & $\varepsilon$ & $(1, 3)$ \\
      $(1, 4, 7)$ & $\varepsilon$ & $ \varepsilon$ & $\varepsilon$ & $(1, 4)$ \\
      \vspace*{0.5cm}
  \end{tabular}
\end{table}

\section{Equi-Distribution of \texorpdfstring{$\luck$ and $\freq_1$}{luck and freq\_1}} \label{secEqui}

In a parking distribution on trees, the drivers that prefer to park at the leaf nodes always park at their preferred nodes. We say a car is \textit{lucky} if it prefers a backbone node and is able to park at its preferred node in the parking process we discuss. Let $\p \in \PK{m}{n}$. We denote the number of lucky drivers of $\p$ by $\luck(\p)$. We also denote the number of cars that prefer to park on the first node (node 1) by $\freq_1(\p)$. These two statistics are often studied in enumeration. We examine these two statistics in $\uu$-parking distributions. Similarly, for $\p' \in \PK{}{n; \uu}$ where $\uu = (1, m+1, 2m+1, \dots)$, we define the two statistics as follows: 
$$\luck(\p') := |i: p'_i = mi - m + 1|$$
$$\freq_j(\p') := |i: p'_i = j|, \  1 \leq j \leq mn - m + 1$$

It is clear that if $\p$ and $\p'$ are isomorphic, i.e., $\theta(\p) = \p'$ then

$$\luck(\p) = \luck(\p')$$
$$\freq_j(\p)  = \begin{cases}
    \freq_j(\p') \text{ if $j\mod m = 1$} \\
    \freq_j(\p') + 1 \text{ otherwise}
\end{cases},  $$

where $\theta$ is the bijection defined in the proof of proposition $\ref{prop:iso}.$
\begin{theorem}
  \label{thm:equi}
  The number of $\uu$-parking distributions $\p$ of length $n$ with $\luck(\p) = k$ is equal to the number of $\uu$-parking distributions $\p'$ of length $n$ with $\freq_1(\p') = k$.
\end{theorem}

\begin{proof}
  We define an involution $\tau: \PK{}{n; \uu} \to \PK{}{n; \uu}$ as follows: Let $\p \in \PK{}{n; \uu}$ and let $(1, \p_1, \dots, \p_{m+1})$ be the first-return decomposition of $\p$. Then define $\tau$ as:
  \[
  \tau(\varepsilon) = \varepsilon
  \]
  \[
  \tau(\p) = (1, \tau(\p_{m+1}), \p_{2}, \p_{3}, \dots, \p_{m}, \tau(\p_{1})).
  \]
  We claim that $\luck(\p) = \freq_{1}(\tau(\p))$ and prove this by induction on $n$. For the base case, when $\p = \varepsilon$, the claim is trivially true. Assume $\luck(\p') = \freq_{1}(\tau(\p'))$ for all $\p' \in \PK{}{l;\uu}$ where $0 \leq l \leq n-1$. Then $\luck(\p) = 1 + \luck(\p_{m+1})$. On the other hand, $\luck(\tau(\p)) = 1 + \luck(\tau(\p_1))$. By the induction hypothesis, $1 + \luck(\p_{m+1}) = 1 + \omega_{1}(\tau(\p_{m+1})) = \freq_{1}(\tau(\p))$ and $1 + \luck(\tau(\p_1)) = 1 + \omega_{1}(\p_{1}) = \freq_{1}(\p)$.
\end{proof}

\newcommand{\R}[1]{R^{(m)}}
We denote the number of $\uu$-parking distributions with length $n$ and $\luck$ value $k$ by $C_{n, k}^{(m)}$. Define a polynomial $\R{m}_n(q)$ as follows:

\[
\R{m}_n(q) = \sum_{\p \in \PK{}{n; \uu}} q^{\luck(\p)} = \sum_{\p \in \PK{}{n; \uu}} q^{\freq_1(\p)} = \sum_{k \geq 0} C_{n, k}^{(m)} q^k.
\]
The polynomial $\R{m}_n(q)$ gives rise to another $q$-analog of the Fuss-Catalan numbers. Table \ref{tab:qanalogs} shows the first few of them for $m = 2$, $3$, and $4$.
\begin{table}
  \centering
\caption{$R_n^{(m)}$}
\label{tab:qanalogs}
  \begin{tabular}{c|lll}
        {$n$}\textbackslash{$m$}&  $2$&  $3$& $4$\\
        \hline
       $0$&  $1$&  $1$& $1$\\
       $1$&  $q$&  $q$& $q$\\
       $2$&  $q^2 + 2q$&  $q^2 + 3q$& $q^2 + 4q$\\
       $3$&  $q^3 + 4q^2 + 7q$&  $q^3 + 6q^2 + 15q$& $q^3 + 8q^2 + 26q$\\
$4$& $q^4 + 6q^3 + 18q^2 + 30q$& $q^4 + 9q^3 + 39q^2 + 91q$&$q^4 + 12q^3 + 68q^2 + 204q$\\
  \end{tabular}  
\end{table}
Let
$$\mathcal{B}_m(x;q) = \sum_{n\geq 0}\R{m}_n(q)x^n,$$
  We use the above generating function in the subsequent sections of the paper. 

\section{Symmetric Joint Distribution of $\luck$ and $\omega_1$}\label{secSymm}
Let $\p \in \PK{}{n; (1, m+1, 2m+1, \dots)}$. Let $f(\p)$ and $g(\p)$ be the indices of $\p$'s first fixed points of type 1 and type $m$ respectively.
\begin{theorem}
\label{bigT}
Define the multivariate polynomial
$$\gamma_{n}^{(m)}(q, t, u, v) = \sum_{\p} q^{\luck(\p)}t^{\omega_1(\p)} u^{f(\p)} v^{g(\p)},$$
where the sum is over all $\p \in \PK{}{n; \uu}$. Let $ [x^{n}]\Gamma_m(x;q, t, u, v) =\gamma_{n}^{(m)}(q, t, u, v)$, then we have

$$\Gamma_m(x;q, t, u, v) = 1 + xq t(uv)^2\mathcal{B}_m^{m-1}(x) \mathcal{B}_m(vx;q)\mathcal{B}_m(uvx;t).$$
    
\end{theorem}
\newcommand{\g}[2]{\gamma^{(#1)}_{#2}}

\begin{proof}
  We decompose $\p$ into $\p_1$, $w$, and $\p_{m+1}$, where $\p_1 \in \PK{}{f(\p)-2; (1, 3, \dots)}$ and $\p_{m+1} \in \PK{}{n - g(\p) + 1; (1, 3, 5, \dots)}$. From this decomposition, we note

 \begin{align*}
 \luck(\p) &= 1 + \luck(\p_{m+1}), \\
 \freq_1(\p) &= 1 + \freq_1(\p_1).
 \end{align*}

Thus, we can rewrite $\g{m}{n}$ as follows:

 \begin{align*}
 \g{m}{n}(q, t, u, v) &= \sum_{k = 2}^{n+1} \sum_{r = k}^{n+1} u^k v^r\sum_{\substack{
\p_1, \p_{m+1}, w\\
|\p_1| = k-2, \\ |\p_{m+1}| = n - r + 1
}} q^{1 + \luck(\p_{m+1})} t^{1 + \freq_1(\p_1)} \\
 &= qt\sum_{k = 2}^{n+1} \sum_{r = k}^{n+1} u^k v^r \sum_{w}\sum_{\p_1} t^{\freq_1(\p_1)}\sum_{\p_{m+1}} q^{\luck(\p_{m+1})},
 \end{align*}

where $\p_1 \in \PK{}{k-2; (1, m+1, \dots)}$ and $\p_{m+1} \in \PK{}{n - r + 1; (1, m+1, \dots)}$, and $w = (p_{k} - k + 1, \dots, p_{r - k} - k + 1)$ is a $\uu$-parking distribution where $\uu = (m-1, 2m-1, \dots)$. By the definition provided in section \ref{secUpar}, $\sum_w 1 = h_{r - k, 1, 1}^{(m)}$, thus we have

 \begin{equation*}
 \g{m}{n}(q, t, u, v) = qt\sum_{k = 2}^{n+1} \sum_{r = k}^{n+1} u^k v^r h_{r - k, 1, 1}^{(m)}\sum_{\p_1} t^{\freq_1(\p_1)}\sum_{\p_{m+1}} q^{\luck(\p_{m+1})}.
 \end{equation*}

We then make the substitution $k \to k - 2$ followed by $r \to r - k + 2$ and rearrange the terms to obtain

 \begin{equation*}
 \g{m}{n}(q, t, u, v) = qt(uv)^2 \sum_{k=0}^{n-1} (uv)^{k} R_{k}^{(m)}(t) \sum_{r=0}^{n-1 - k} h^{(m)}_{r, 1, 1} v^{r}  R_{n -1 -k-r}^{(m)}(q),
 \end{equation*}
where the substitution of $\R{m}_{k}$ and $\R{m}_{n-1-k-r}$ follows from theorem \ref{thm:equi}. We can apply Cauchy's product rule and theorem \ref{thm:rec} to get

 \begin{equation*}
 \Gamma_m(x;q, t, u, v) = 1 + xqt(uv)^2\mathcal{B}_m^{m-1}(x) \mathcal{B}_m(vx;q)\mathcal{B}_m(uvx;t).
 \end{equation*}
\end{proof}

\newtheorem{cor}{Corollary}
\begin{cor}
  \label{corinv}
    $$\mathcal{B}_m(x;q) = \sum_{n \geq 0} \sum_{\p} q^{\luck(\p)}x^n=  \dfrac{1}{1 - qx \mathcal{B}_m(x)}$$
\end{cor}
\begin{proof}
    Substitute $u = v = 1$ in the functional equation of theorem \ref{bigT}.
    
\end{proof}

\begin{cor}
  \label{corsym}
    If $\p$ is a parking distribution on a $m$-regular caterpillar tree $T$, then the number of lucky drivers and the number of drivers that prefer to park at node $1$ form a symmetric joint distribution. In fact, the polynomial $\g{m}{n}(q, t, 1, 1)/qt$ is a linear combination of complete homogeneous polynomials on $q$ and $t.$
    
\end{cor}

\begin{proof}
  We observe that 

\[
\begin{array}{rcl}
\dfrac{q \mathcal{B}_m(x; q) - t \mathcal{B}_m(x; t)}{q - t} & = & \dfrac{q \mathcal{B}_m(x; q) - t \mathcal{B}_m(x; t)}{q - t} \cdot \dfrac{\mathcal{B}_m(x;q)\mathcal{B}_m(x;t) }{\mathcal{B}_m(x;q)\mathcal{B}_m(x;t) } \\
& = & \dfrac{\mathcal{B}_m(x;q)\mathcal{B}_m(x;t)}{q - t} \bigg(\dfrac{q}{\mathcal{B}_m(x;t)} - \dfrac{t}{\mathcal{B}_m(x;q)}\bigg)\\
& = & \dfrac{\mathcal{B}_m(x;q)\mathcal{B}_m(x;t)}{q - t} \cdot (q - t) \\
& = & \mathcal{B}_m(x;q)\mathcal{B}_m(x;t), \\
\end{array}
\]
where we used corollary \ref{corinv} to reduce $\frac{q}{\mathcal{B}_m(x;t)} - \frac{t}{\mathcal{B}_m(x;q)}$ to $q - t$. This implies

\[
\sum_{k = 0}^{n} \R{m}_k(q)\R{m}_{n-k}(t) = \sum_{k = 1}^{n}C_{n, k}^{(m)}\dfrac{q^{k+1} - t^{k+1}}{q - t} = \sum_{k = 1}^n C_{n, k}^{(m)} \sum_{s = 0}^k q^{k-s}t^{s}
\]

Consecutively,

\[
\g{m}{n}(q, t, 1, 1) = qt \sum_{r=0}^{n-1} h^{(m)}_{r, 1, 1} \sum_{k=0}^{n-1 - r} R_{k}^{(m)}(t) R_{n -1 -k-r}^{(m)}(q)
\]

\[
= qt \sum_{r=0}^{n-1} h^{(m)}_{r, 1, 1}\sum_{k = 0}^{n-1-r}C_{n-r-1, k}^{(m)} \sum_{s = 0}^k q^{k-s}t^{s}
\]

\[
= qt \sum_{r=0}^{n-1} h^{(m)}_{r, 1, 1}\sum_{k = 0}^{n-1-r}C_{n-r-1, k}^{(m)} h_k(q, t),
\]

where \(h_k\) is the complete symmetric homogeneous polynomial of degree \(k\).

\end{proof}

Before we proceed to the next section, we would like show how $\gamma_n^{(m)}$ looks like for a few cases when $u = v = 1$.

\begin{table}[h!]
  \centering
  \caption{$\gamma_n^{(2)}$}
  \label{tab:g2}
  \begin{tabular}{c|p{6cm}|l}
      $n$ & $\gamma_n^{(2)}(qt)/qt$ &   \\
      \hline
      $1$ & $1$ & $h_0(q, t)$ \\
      $2$ & $t + q + 1$ & $h_1(q, t) + h_0(q, t)$ \\
      $3$ & $t^2 + qt + q^2 + 3t + 3q + 3$& $h_2(q, t) + 3h_1(q, t) + 3h_0(q, t)$\\
      $4$ & $t^3 + t^2q + tq^2 + q^3 + 5t^2 + 5qt + 5q^2 + 12t + 12q + 12$ & $h_3(q, t) + 5h_2(q, t) + 12h_1(q, t) + 12h_0(q, t)$\\
  \end{tabular}
\end{table}

\begin{table}[h!]
  \centering
  \caption{$\gamma_n^{(3)}$}
  \label{tab:g3}
  \begin{tabular}{c|p{6cm}|l}
      $n$ & $\gamma_n^{(3)}/qt$  &\\
      \hline
      $1$ & $1$  &$h_0(q, t)$\\
      $2$ & $t + q + 2$ &$h_1(q, t) + 2h_0(q, t)$ \\
      $3$ & $t^2 + tq + q^2 + 5t + 5q + 9$  &$h_2(q, t)  + 5h_1(q, t) + 9h_0(q, t)$\\
      $4$ & $t^3 + t^2q + tq^2 + q^3 + 8t^2 + 8tq + 8q^2 + 30t + 30q + 52$  &$h_3(q, t) + 8h_2(q, t) + 30h_1(q, t) + 52h_0(q, t)$\\
  \end{tabular}
\end{table}

\begin{table}[h!]
  \centering
  \caption{$\gamma_n^{(4)}$}
  \label{tab:g4}
  \begin{tabular}{c|p{6cm}|l}
      $n$ & $\gamma_n^{(4)}/qt$  &\\
      \hline
      $1$ & $1$  &$h_0(q, t)$\\
      $2$ & $t + q + 3$ &$h_1(q, t) + 3 h_0(q, t)$\\
      $3$ & $t^2 + qt + q^2 + 7t + 7q + 18$  &$h_2(q, t)  + 7h_1(q, t) + 18h_0(q, t)$\\
      $4$ &$t^3 + t^2q + tq^2 + q^3 + 11t^2 + 11qt + 11q^2 + 56t + 56q + 136$  &$h_3(q, t) + 11h_2(q, t) + 56h_1(q, t) + 136h_0(q, t)$\\
  \end{tabular}
\end{table}

\newpage
\newcommand{\Si}{\mathcal{S}}

\section{An Abstraction} \label{secGen}

The proof of theorem \ref{bigT} suggests that if we have $m+1$ statistics $\Si_0, \Si_1, \dots, \Si_{m}$, each of which is equidistributed with $\luck$, and if each $\Si_i$ satisfies $\Si_i(\p) = \Si_i(\p_{i+1}) + C_i$ for some constant $C_i$, where $\ip_{i+1}$ is as defined in section \ref{secFTD}, then the symmetric joint distribution of these $m+1$ statistics holds.

\begin{example}
  \label{exam}
  As we shall prove in this section, $\omega_1 \sim \omega_2 \sim \dots \sim \omega_m \sim \luck$. When we look at the terms of the polynomial $\sum_\p q_0^{\luck(\p)}q_1^{\omega_1(\p)} q_2^{\omega_2(\p)}$ where $\p = \PK{2}{n}$ for different values of $n$ we find the following.

  \begin{table}[h!]
      \centering
      \caption{Polynomial expansions for different values of $n$}
      \label{tab:g4}
      \begin{tabular}{c|l}
          $n$ & $\sum_\p q_0^{\luck(\p)}q_1^{\omega_1(\p)} q_2^{\omega_2(\p)}$, $\p = \PK{2}{n}$ \\
          \hline
          $1$ & $q_0q_1$\\
          $2$ &$q_0q_1q_2[h_1(q_0, q_1, q_2)]$\\
          $3$ &$q_0q_1q_2[h_2(q_0, q_1, q_2) + 2h_1(q_0, q_1, q_2)]$\\
          $4$ &$q_0q_1q_2[h_3(q_0, q_1, q_2) + 4h_2(q_0, q_1, q_2) + 7h_1(q_0, q_1, q_2)]$\\
      \end{tabular}
  \end{table}
\end{example}

where $h_n$ is the complete homogeneous polynomial of degree $n$. Compare the above with $R_n^{(2)}$. Using this observation, we generalize our first result as follows.

\begin{theorem}
  Let $\p \in \PK{}{n; (1, m+1, 2m+1, \dots)}$, and let $\p = (1, \p_1, \dots, \p_{m+1})$ be the first-return decomposition of $\p$. Let $\Si_0, \dots, \Si_{m}$ be $m+1$ statistics on the set $\bigcup_{i \geq 0} \PK{}{i; (1, m+1, 2m+1, \dots)}$, where each $\Si_i: \PK{}{n; (1, m+1, 2m+1, \dots)} \mapsto \mathbb{N}$ is equidistributed with the $\luck$ statistic, i.e.,
  $$\Si_i \sim \luck.$$
  If $\Si_i(\p) = \Si_i(\p_{i+1}) + C_i$ for all $n \geq 1$ and some constants $C_i$ depending only on $i$, and if
  $$\Gamma_m(x; q_0, \dots, q_{m}) = \sum_{n \geq 0} x^n \sum_{\p} q_0^{\Si_0(\p)} \cdots q_{m}^{\Si_{m}(\p)},$$
  then we have
  $$\Gamma_m(x; q_0, \dots, q_{m})  = 1 + x \prod_{i = 0}^{m+1} q_i^{C_i}\mathcal{B}_m(x;q_i).$$
\end{theorem}

\begin{proof}
  The proof is similar to that of theorem \ref{bigT}. If $\gamma^{(m)}_{n}(q_0, \dots, q_m) = [x^n] \Gamma_m(x; q_0, \dots, q_{m})$, then

\[
\begin{array}{rcl}
\gamma^{(m)}_{n}(q_0, \dots, q_m) & = & \displaystyle\sum_{\p_1} \cdots \displaystyle\sum_{\p_{m+1}} q_0^{\Si_0(\p)} \cdots q_{m}^{\Si_{m}(\p)} \\
& = & \displaystyle\sum_{\p_1} \cdots \displaystyle\sum_{\p_{m+1}} q_0^{\Si_0(\p_1) + C_0} \cdots q_{m}^{\Si_{m}(\p_{m+1}) + C_m} \\
& = & \displaystyle\prod_{i=0}^{m} \sum_{\p_i} q_i^{\Si_i(\p_i) + C_i} \\
& = & \displaystyle\prod_{i = 0}^{m} q_{i}^{C_i} \sum_{\p_i} q_i^{\Si_i(\p_i)} \\
& = & \displaystyle\prod_{i = 0}^{m} q_{i}^{C_i} \sum_{\alpha} \R{m}_{\alpha_i}(q_i),
\end{array}
\]
where the sum is over all weak $(m+1)$-compositions $\alpha = (\alpha_0, \dots, \alpha_m)$ of $n-1$. The theorem follows after applying Cauchy's product rule.
\end{proof}

\begin{cor}
  Let $\mathbf{q} \in \PK{m}{n}$, $\mathcal{R}_0(\mathbf{q}) = \luck(\p)$, and $\mathcal{R}_i(\mathbf{q}) = \freq_{i}(\mathbf{q})$ for $1 \leq i \leq m$. Then we have
  $$\sum_{n \geq 0} x^n \sum_{\mathbf{q} \in \PK{m}{n}} \prod_{i = 0}^{m}q_{i}^{\mathcal{R}_i(\mathbf{q})} = 1 + x \prod_{i = 0}^{m+1} q_i\mathcal{B}_m(x;q_i).$$
\end{cor}

\begin{proof}
  Let $\p \in \PK{}{n; (1, m+1, 2m+1, \dots)}$. The left-hand side of the equation is equal to $\Gamma_m(x; q_0, \dots, q_m)$ if $\mathcal{S}_0(\p) = \mathcal{R}_0(\theta(\p))$, $\mathcal{S}_1(\p) = \mathcal{R}_1(\theta(\p))$, and $\mathcal{S}_i(\p)+1 = \mathcal{R}_i(\theta(\p))$ for $2 \leq i \leq m$, with $\theta$ defined in the proof of Proposition \ref{prop:iso}.

  Given the first-return decomposition of $\p$, define a mapping $\eta$ as follows: Let $\p'$ be a parking distribution initially containing only $1$. Insert $\freq_1(\p_1)$ $1$s, $\freq_1(\p_2)$ $2$s, and so on until all $\freq_1(\p_{m+1})$ $m+1$s are inserted in $\p'$. Then, starting from the right-most non-empty $\p_j$, for each parking preference $p_{i, j}$ of $\p_j$ different from $1$, insert $p_{i, j} + m(1 + \sum_{k>j} |\p_k|)$ into $\p'$. Define the resulting parking distribution as $\eta(\p)$. Table \ref{tab:eta} shows the values of $\eta$ for $\PK{}{3; (1,3, \dots)}$.

\begin{table}[h!]
  \centering
  \caption{Example values of $\eta(\p)$ for $\p \in \PK{}{3; (1,3, \dots)}$.}
  \label{tab:eta}
  \begin{tabular}{c|ccc|c}
  $\p$ & $\p_1$ & $\p_2$ & $\p_3$ & $\eta(\p)$\\
  \hline
  \hline
  $(1, 1, 1)$ & $(1, 1)$ & $\varepsilon$ & $\varepsilon$ & $(1, 1, 1)$ \\
  $(1, 1, 2)$ & $(1, 2)$ & $\varepsilon$ & $\varepsilon$ & $(1, 1, 4)$ \\
  $(1, 1, 3)$ & $(1, 3)$ & $\varepsilon$ & $\varepsilon$ & $(1, 1, 5)$ \\
  $(1, 1, 4)$ & $(1)$   & $(1)$ & $\varepsilon$ & $(1, 1, 2)$ \\
  $(1, 1, 5)$ & $(1)$   & $\varepsilon$ & $(1)$ & $(1, 1, 3)$ \\
  $(1, 2, 2)$ & $\varepsilon$     & $(1, 1)$ & $\varepsilon$ & $(1, 2, 2)$ \\
  $(1, 2, 3)$ & $\varepsilon$     & $(1, 2)$ & $\varepsilon$ & $(1, 2, 4)$ \\
  $(1, 2, 4)$ & $\varepsilon$     & $(1, 3)$ & $\varepsilon$ & $(1, 2, 5)$ \\
  $(1, 2, 5)$ & $\varepsilon$     & $(1)$ & $(1)$ & $(1, 2, 3)$ \\
  $(1, 3, 3)$ & $\varepsilon$     & $\varepsilon$ & $(1, 1)$ & $(1, 3, 3)$ \\
  $(1, 3, 4)$ & $\varepsilon$     & $\varepsilon$ & $(1, 2)$ & $(1, 3, 4)$ \\
  $(1, 3, 5)$ & $\varepsilon$     & $\varepsilon$ & $(1, 3)$ & $(1, 3, 5)$
  \end{tabular}
\end{table}

This mapping is a bijection because we can uniquely reconstruct $\p$ from $\eta(\p)$ by first reconstructing the first-return decomposition of $\p$: Place $\omega_j(\p)$ $1$s in $\p_{j}$ for $1 \leq j \leq m + 1$. For the remaining elements $p'_i$ of $\eta(\p)$, insert $q_{i, j} = p'_i - m(1 + \sum_{k>j} |\p_k|)$ into the largest $j$ such that $q_{i, j} > 0$ and $\p_j$ is a valid $\uu$-parking distribution after inserting $q_{i, j}$.

This proves that for $1 \leq i \leq m+1$ and $\p \in \PK{}{n;(1, m+1, \dots)}$, $\omega_{k}(\p) \sim \omega_1(\p_k)$. As a result, $\omega_k \sim \luck$ and is enumerated by $\omega_1(\p_k)$. This completes the proof of the corollary.

\end{proof}

As with Corollary \ref{corsym}, we can express $\gamma_n^{(m)}(q_0, \dots, q_{m})$ as a linear combination of complete homogeneous polynomials in $q_0, \dots, q_{m}$. As this is central to our paper, we present it as a theorem.

\begin{theorem}
    Let $\Si_0, \dots, \Si_m$ be $m+1$ statistics on $\PK{}{n; (1, m+1, 2m+1, \dots)}$ that are equidistributed with the $\luck$ statistic, i.e., $\Si_i \sim \luck$, and assume that $\Si_i(\p) = \Si_i(\p_i) + C_i$ for each statistic. Then the polynomial
    $$\gamma_n^{(m)}(q_0, \dots, q_m) = \sum_{\p \in \PK{m}{n}} \prod_{i = 0}^m q_i^{\Si_i(\p)}$$
    is a linear combination of complete homogeneous polynomials in $q_0, \dots, q_{m}$.
\end{theorem}

\begin{proof}
  The proof of Corollary \ref{corsym} establishes that

  \begin{equation}
      \sum_{k = 0}^{n} \R{m}_k(q_0)\R{m}_{n-k}(q_1) = \sum_{k = 1}^n C_{n, k}^{(m)} \sum_{s = 0}^k q_0^{k-s}q_1^{s}. \label{eq2}
  \end{equation}

  We use induction to show that

  \begin{equation}
      P_{n, t}(q_0, \dots, q_{t-1}) = \sum_{\alpha_n} \R{m}_{\alpha_{n, 1}}(q_0)\cdots\R{m}_{\alpha_{n, t}}(q_{t-1}) = \sum_{k = 1}^n C_{n, k}^{(m)} \sum_{\alpha_k} \prod_{i = 0}^t q_{i}^{\alpha_{k, i}} \label{eq1}
  \end{equation}
  where $\alpha_n = (\alpha_{n, 1}, \dots, \alpha_{n, t})$ is a $t$-weak-composition of $n$ for $2 \leq t \leq m+1$. The base case for $t = 2$ is proved by Corollary \ref{corsym}.

  For the induction step, assume that \eqref{eq1} is true for $t = k$:
  \[
  P_{n, k} = \sum_{\alpha_n} \R{m}_{\alpha_{n, 1}}(q_0) \cdots \R{m}_{\alpha_{n, k}}(q_k) = \sum_{i = 0}^{n} P_{i, k} \cdot \R{m}_{n-i}(q_{k}).
  \]
  Then,

  \[
  P_{n, k+1} = \sum_{i = 0}^{n} \left( \sum_{r = 1}^i C_{i, r}^{(m)} \sum_{\alpha_r} \prod_{j = 0}^{k-1} q_{j}^{\alpha_{r, j+1}} \right) \cdot \left( \sum_{s = 0}^{n-i} C_{n-i, s}^{(m)} q_{k}^s \right).
  \]

  Let $\mathcal{Q}(r, k-1) = \sum_{\alpha_r} \prod_{j = 0}^{k-1} q_{j}^{\alpha_{r, j+1}}$. Then \eqref{eq1} is equivalent to:

  \[
  \sum_{k = 0}^n \sum_{i = 0}^k C_{k, i}^{(m)} \mathcal{Q}(i, k-1) \cdot \sum_{s = 0}^{n-i} C_{n-i, s}^{(m)} q_k^s.
  \]

  We can simplify this to:

  \[
  \sum_{k = 0}^n \sum_{i = 0}^k C_{k, i}^{(m)} \mathcal{Q}(i, k-1) \cdot \sum_{s = 0}^{n-i} C_{n-i, s}^{(m)} q_k^s = \sum_{k = 0}^n C_{n, k}^{(m)} \mathcal{Q}(k, k).
  \]

  Each term on the right-hand side results from a product of one term from each factor on the left-hand side. This confirms:

  \[
  \sum_{k = 0}^n \sum_{i = 0}^k C_{k, i}^{(m)} \mathcal{Q}(i, s) \cdot \sum_{s = 0}^{n-i} C_{n-i, s}^{(m)} q_s^k = \sum_{k = 0}^n C_{n, k}^{(m)} \mathcal{Q}(k, s+1).
  \]

  The theorem follows by taking $s = k$ and applying the induction hypothesis.
\end{proof}

If $B_m(n, k_0, \dots, k_m)$ is the number of parking distributions on $\Cat{m}{n}$ with $k_0$ lucky cars and $k_j$ preferences of node $j$ for $1 \leq j \leq m$ then we proved that

$$B_m(n, k_0, \dots, k_m) = B_m(n, l_0, \dots, l_m)$$
if $\sum_i {l_i} = \sum_i k_i$. For instance, $B_2(4, k_0, k_1, k_2)$ is shown in table \ref{tab:tensor}

\begin{table}[h!]
  \centering
  \caption{$B_2(4, k_0, k_1, k_2)$}
  \vspace{.3cm}
  \label{tab:tensor}
  \begin{tabular}{c|c|cccc}
      & & \multicolumn{4}{c}{$k_2$} \\
      \cline{3-6}
      $k_0$ & $k_1$ & 1 & 2 & 3 & 4 \\
      \hline
      \multirow{4}{*}{1}
      & 1 & 0 & 7 & 4 & 1 \\
      & 2 & 7 & 4 & 1 & 0 \\
      & 3 & 4 & 1 & 0 & 0 \\
      & 4 & 1 & 0 & 0 & 0 \\
      \hline
      \multirow{4}{*}{2}
      & 1 & 7 & 4 & 1 & 0 \\
      & 2 & 4 & 1 & 0 & 0 \\
      & 3 & 1 & 0 & 0 & 0 \\
      & 4 & 0 & 0 & 0 & 0 \\
      \hline
      \multirow{4}{*}{3}
      & 1 & 4 & 1 & 0 & 0 \\
      & 2 & 1 & 0 & 0 & 0 \\
      & 3 & 0 & 0 & 0 & 0 \\
      & 4 & 0 & 0 & 0 & 0 \\
      \hline
      \multirow{4}{*}{4}
      & 1 & 1 & 0 & 0 & 0 \\
      & 2 & 0 & 0 & 0 & 0 \\
      & 3 & 0 & 0 & 0 & 0 \\
      & 4 & 0 & 0 & 0 & 0 \\
  \end{tabular}
\end{table}

\section{Conclusion}
\label{sec:conclusion}
In this work, we have presented an analysis of the symmetric joint distribution of the $\luck$ and $\omega_1$ statistics in caterpillar parking distributions. We have shown that this symmetry can be generalized for any $m+1$ statistics that are equidistributed with the $\luck$ statistic and depend only on the value of the statistic for the $i$-th part in the first-return decomposition of the parking distribution. We have also demonstrated that the $q$-analog of these statistics can be expressed using complete symmetric homogeneous polynomials.

Furthermore, we have derived a $q, t$-analog for the Fuss-Catalan numbers and established a linear combination of complete homogeneous polynomials in $q_0, \dots, q_m$ as a closed form expression for $\gamma_n^{(m)}(q_0, \dots, q_m)$.

We conclude this paper by stating some unsolved problems.
\begin{enumerate}
  \item Is there a closed formula for the coefficient of $q^at^bv^cu^d$ in $\gamma_n^{(m)}(q,t,v,u)$? What about the coefficient pf $\prod_i q_i^{a_i}$ in $\gamma_n^{(m)}(q_0, \dots, q_m)$?
  \item Can we express $\gamma_n^{(m)}$ in terms of elementary/power/Schur symmetric functions?
  \item Is there a closed expression for the $q, t$-analogs discussed in this paper?
\end{enumerate}

\bibliographystyle{amsplain}
\bibliography{sample}

\providecommand{\bysame}{\leavevmode\hbox to3em{\hrulefill}\thinspace}
\providecommand{\MR}{\relax\ifhmode\unskip\space\fi MR }
\providecommand{\MRhref}[2]{%
  \href{http://www.ams.org/mathscinet-getitem?mr=#1}{#2}
}
\providecommand{\href}[2]{#2}
\begin{thebibliography}{10}

\bibitem{Multivar}
Jean-Christophe Aval, \emph{Multivariate {F}uss--{C}atalan numbers}, Discrete
  Mathematics \textbf{308} (2008), no.~20, 4660--4669.

\bibitem{flattendcatalan}
Jean-Luc Baril, Pamela~E. Harris, and Jos{\'e}~L. Ram{\'i}rez, \emph{Flattened
  {C}atalan words}, arXiv preprint arXiv:2405.05357, 2024, Available at
  \url{https://arxiv.org/abs/2405.05357}.

\bibitem{Gpark}
Brian Benson, Deeparnab Chakrabarty, and Prasad Tetali, \emph{$g$-parking
  functions, acyclic orientations and spanning trees}, Discrete Mathematics
  \textbf{310} (2010), no.~8, 1340--1353.

\bibitem{Butler_2017}
Steve Butler, Ron Graham, and Catherine~H. Yan, \emph{Parking distributions on
  trees}, European Journal of Combinatorics \textbf{65} (2017), 168--185.

\bibitem{hypergraph}
Parth Chavan, Andrew Lee, and Karthik Seetharaman, \emph{Hypergraph
  {F}uss-{C}atalan numbers}, arXiv preprint arXiv:2202.01111, 2022, Available
  at \url{https://arxiv.org/abs/2202.01111}.

\bibitem{Deutsch_1999}
Emeric Deutsch, \emph{An involution on {D}yck paths and its consequences},
  Discrete Mathematics \textbf{204} (1999), no.~1--3, 163--166.

\bibitem{Franon1975}
Jean Fran\c{c}on, \emph{Acyclic and parking functions}, Journal of
  Combinatorial Theory, Series A \textbf{18} (1975), no.~1, 27--35.

\bibitem{Garsia_2001}
A.~M. Garsia and J.~Haglund, \emph{A positivity result in the theory of
  {M}acdonald polynomials}, Proceedings of the National Academy of Sciences
  \textbf{98} (2001), no.~8, 4313--4316.

\bibitem{Garsia1996ARQ}
Adriano~M. Garsia and Mark Haiman, \emph{A remarkable $q$, $t$-{C}atalan
  sequence and $q$-{L}agrange inversion}, Journal of Algebraic Combinatorics
  \textbf{5} (1996), 191--244.

\bibitem{luckycars}
Ira~M. Gessel and Seunghyun Seo, \emph{A refinement of {C}ayley's formula for
  trees}, arXiv preprint arXiv:math/0507497, 2005, Available at
  \url{https://arxiv.org/abs/math/0507497}.

\bibitem{Haiman_2002}
Mark Haiman, \emph{Vanishing theorems and character formulas for the {H}ilbert
  scheme of points in the plane}, Inventiones Mathematicae \textbf{149} (2002),
  no.~2, 371--407.

\bibitem{Konheim_1966}
Alan~G. Konheim and Benjamin Weiss, \emph{An occupancy discipline and
  applications}, SIAM Journal on Applied Mathematics \textbf{14} (1966), no.~6,
  1266--1274.

\bibitem{Kung2003}
Joseph P.~S. Kung and Catherine Yan, \emph{{G}onarov polynomials and parking
  functions}, Journal of Combinatorial Theory, Series A \textbf{102} (2003),
  no.~1, 16--37.

\bibitem{density}
Wojciech Mlotkowski, Karol~A. Penson, and Karol Zyczkowski, \emph{Densities of
  the {R}aney distributions}, Documenta Mathematica \textbf{18} (2013),
  1573--1596.

\bibitem{hopfalg}
Jean-Christophe Novelli and Jean-Yves Thibon, \emph{A {H}opf algebra of parking
  functions}, arXiv preprint arXiv:math/0312126, 2003, Available at
  \url{https://arxiv.org/abs/math/0312126}.

\bibitem{StanNonCrossing}
Richard~P. Stanley, \emph{Parking functions and noncrossing partitions},
  Electronic Journal of Combinatorics \textbf{4} (1996), no.~2, R20.

\bibitem{Stanley1998}
\bysame, \emph{Hyperplane arrangements, parking functions and tree inversions},
  pp.~359--375, Birkh\"{a}user Boston, 1998.

\bibitem{YanGeneral}
Catherine~H. Yan, \emph{Generalized parking functions, tree inversions, and
  multicolored graphs}, Advances in Applied Mathematics \textbf{27} (2001),
  no.~2--3, 641--670.

\end{thebibliography}

\end{document}